\magnification=1200 \baselineskip=12pt \font\abf=cmbx12
\font\arm=cmr12 \centerline{\abf  Existence of Gradient
K\"ahler-Ricci Solitons}
\bigskip
\centerline {{\arm Huai-Dong Cao}\footnote* {Research partially
supported by a NSF grant and  a Sloan Research Fellowship.}}

\vskip 0.5in \noindent {\S1. \bf Introduction}
\bigskip
We consider the K\"ahler-Ricci flow
$${\partial \over {\partial t}}g_{i\bar j} =- R_{i\bar j}\eqno (1.1)$$
on a complete noncompact  $n$-dimensional K\"ahler manifold $M^n$.
A solution $g_{i\bar j}(t)$, for $t\geq 0$, is called a
K\"ahler-Ricci soliton if it moves along Eq.(1.1) under a
one-parameter family of biholomorphisms. If this comes from a
holomorphic vector field $V=(V^i)$ then we have a soliton when
$$V_{i,\bar j}+ V_{\bar j,i}=R_{i\bar j},$$
since the metric changes by its Lie derivatives. When the
holomorphic vector field is the gradient of a function we say we
have a gradient soliton. In this case, we have
$$  R_{i\bar j}=f_{,i\bar j} \qquad \hbox {and} \qquad  f_{,ij}=0 \eqno (1.2)$$
for some real-valued  function $f$ on $M$. We remark that the
condition $ f_{,ij}=0 $ is equivalent to saying that the vector
field
 $V=\bigtriangledown f$ is
holomorphic.

On a compact K\"ahler manifold $M^n$  with positive first Chern
class $C_1(M)$, we can consider the normalized K\"ahler-Ricci flow
$${\partial \over {\partial t}}g_{i\bar j} =- R_{i\bar j}+g_{i\bar
j}\eqno (1.3)$$ A gradient K\"ahler-Ricci soliton of Eq.(1.3) is
then characterized by the equations
 $$R_{i\bar j}-g_{i\bar j}=f_{,i\bar j} \qquad {\hbox {and}}\qquad f_{,ij}=0. \eqno (1.4)$$

The notion of  Ricci soliton was first introduced by Hamilton [4].
It has played a very important role in the study of Ricci flow.
One of the important aspects is its relation with the Harnack
estimates (see  [3] and [5]). Another one is the fact that  Ricci
solitons  often arise as the limits of singularities in the Ricci
flow (see [6]). However, few examples were known. In the
noncompact case, Hamilton [4] wrote down the first example of a
Ricci soliton, called the cigar soliton, on the complex plane
${\bf C}$. It has the form
$$ds^2={|dz|^2\over {1+|z|^2}},\qquad z\in {\bf C} $$
which flows towards the origin by conformal dilations and has
positive Gaussian curvature. It is asymptotic to a flat cylinder
at infinity and has maximal Gaussian curvature at the origin.
Later, Robert Bryant [1] found a gradient Ricci soliton on ${\bf
R}^3$ with positive curvature operator.

Our  purpose in this note is to provide  examples of
K\"ahler-Ricci solitons. 

\proclaim  Theorem 1. For each $n\ge 1$, there exists on ${\bf
C}^n$
  a complete rotationally symmetric gradient K\"ahler-Ricci
soliton of positive sectional curvature. Moreover such a soliton
is unique up to scaling and dilation.

\medskip
\noindent {\bf Remark 1}: For $n=1$, the soliton is just the cigar
soliton found by Hamilton.
\smallskip
 \noindent {\bf Remark 2}: It can be shown that the soliton metric $g$
on {\bf $C^n$} ($n\geq 2$) satisfies the following properties:\
(i) Let $\rho$ denote the distance function from the origin with
respect to $g$, then the volume of the geodesic ball $B_g(0,\rho)$
with respect to the metric $g$ grows like ${\rho}^n$; \ (ii) The
scalar curvature $R$ of the metric $g$ decays like $1/ \rho$.

One interesting open problem is the following uniqueness question:
\medskip

 \noindent {\bf Question}: Suppose $g$ is a complete gradient K\"ahler-Ricci soliton of positive curvature on ${\bf C}^n$, is it true that $g$ is necessarilly a rotationally symmetric one (hence given by Theorem 1)?

The anwser to the above question is affirmative when $n=1$. This
is because in general we always have a Killing vector field given
by $JV$ which gives a $S^1$ action on the soliton. Here $J$
denotes the complex structure of ${\bf C}^n$.

\proclaim  Theorem 2. The total space $X_n$ of the anticanonical
(or canonical) line bundle of the complex projective space ${\bf
P}^{n-1}$ ($n\geq 2$) admits
  a complete rotationally symmetric gradient K\"ahler-Ricci
soliton.  

\medskip
\noindent {\bf Remark 3}: For $n=2$, the manifold $X_2$ is the
tangent bundle to the sphere $S^2={\bf P}^1$ on which Eguch and Hanson [10] constructed the well-known Hyper-K\"ahler (Calabi-Yau) metric. 
The soliton metric restricted to each fiber is quasi-isometric to the cigar soliton
on the complex plane ${\bf C}$. Moreover, Calabi [9] constructed Hyper-K\"ahler metrics on  $X_n$ for $n\ge 3$. 
\medskip

In the compact case, it is known that on Riemann surfaces there
are no Ricci solitons except those of constant curvature.  Tom
Ivey [7] showed that there are no Ricci solitons other than
constant curvature metrics on  compact $3$-dimensional Riemannian
manifolds. This  naturally leads to the question of the existence
of nontrivial Ricci solitons on compact manifolds. On a compact
K\"ahler manifold with positive first Chern class, the
nonvanishing of the Futaki invariant  is a necessary condition for
the existence of a gradient K\"ahler-Ricci soliton which is not
K\"aler-Einstein.  In this article we also construct gradient
K\"ahler-Ricci solitons on certain compact K\"ahler manifolds:

\proclaim  Theorem 3. There exists a unique $U(n)$-invariant (shrinking) 
K\"ahler-Ricci soliton on the total
space $M_k$ of the projective line bundle ${\bf P}(L^k\bigoplus
L^{-k})\buildrel \pi \over \longrightarrow {\bf P^{n-1}}$ for each positive integer $1\le k\le n-1$.

\medskip
\noindent After we finished the construction at Columbia University, we learned from Professor S.
Bando, while he was visiting Courant Institute in 1991, that the existence
of gradient solitons on these compact manifolds had already been
shown by Koiso [8].  Our treatment here, however, is more
explicit. Moreover, we show that the soliton metric on $M_1$ has positive Ricci curvature. 

\medskip
In general, the Ricci soliton equations (1.2) and (1.4) are
nonlinear systems and very difficult to solve. In fact, all known
examples of soliton metrics so far are rotationally symmetric ones
and  both (1.2) and (1.4) are reduced to certain nonlinear ODEs.
The gradient solitons in Theorem 1 and Theorem 2 are  found by
solving these ODEs.


We are very grateful to  Professor Richard Hamilton for bringing
these problems to our attention. We would like to thank Professor
Gang Tian for very helpful conversations and for his interest. We
also like to thank the Institute for Advanced Study for its
support, where part of the work was carried out. Finally, we want
to thank the referee for his careful reading and helpful
suggestions.

\vskip 0.3in

\noindent {\bf \S 2. The Gradient Soliton on ${\bf C}^n$}
\bigskip
In this section we consider rotationally symmetric  K\"akhler
metrics on $\bf C^n$ and derive the Ricci soliton equation for
such metrics.

Any K\"ahler metric $g_{i\bar j}dz^idz^{\bar j}$ on $\bf C^n$
invariant under the unitary group $U(n)$ can be gererated by a
K\"ahler potential $\Phi (z,\bar z)$:
$$g_{i\bar j}={\partial}_i{\partial}_{\bar j}\Phi (z,\bar z)  \eqno (2.1) $$
such that   $$ \Phi (z,\bar z) =w(|z|^2),\qquad  z\in {\bf C^n}
$$ for some smooth function $w(r)$.

In fact it is more convenient if we use the variable $t=\log
|z|^2$. Then we can express the K\"ahler potential as
$$ \Phi (z,\bar z) =u(t),\qquad \bigl(t=\log |z|^2\bigr) \eqno (2.2)$$
 where $u(t)$ is a smooth function  on $(-\infty, \infty)$ and as
 $t\rightarrow -\infty$, it has an expansion
$$u(t)=a_0+a_1e^t+a_2e^{2t}+\cdots, \qquad a_1>0. \eqno (2.3)$$

Conversely, any smooth function $u(t)$ on $(-\infty, \infty)$
 generates a $U(n)$-invariant K\"ahler metric on $\bf C^n$
if and only if it satisfies the differential inequalities
$$ u'(t)>0,\quad u''(t)>0, \qquad t\in (-\infty, \infty) \eqno (2.4)$$
and the asymptotic condition (2.3).

By (2.1) and (2.2) we  have, for $t=\log |z|^2$,
$$g_{i\bar j}={\partial}_i{\partial}_{\bar j}u(t)=e^{-t}u'(t){\delta}_{i\bar j}
+e^{-2t}{\bar z}_iz_j\bigl(u''(t)-u'(t)\bigr). \eqno (2.5)$$

Consequently,
$$g^{i\bar j}=e^t(u'(t))^{-1}{\delta}^{i\bar j}
+ z_i{\bar z}_j\bigl(u''(t)^{-1}-u'(t)^{-1}\bigr) \eqno (2.6)$$
and
$$ \det (g_{i\bar j})=e^{-nt}(u'(t))^{n-1}u''(t) \eqno (2.7)$$

Let
$$f(t)=-\log \det (g_{i\bar j})=nt-(n-1)\log u'(t)-\log u''(t) \eqno (2.8)$$
Then the Ricci tensor of the metric $g_{i\bar j} $ is
$$R_{i\bar j}={\partial}_i{\partial}_{\bar j}f(t)\eqno (2.9)$$

Consider the gradient vector field $V^i=g^{i\bar j}f_{,\bar j}$.
Using (2.6) we get
$$V^i=g^{i\bar j}e^{-t}z_jf'(t)=z_i{f'(t)\over u''(t)}\eqno (2.10)$$
Since $z_i$ is holomorphic and ${f'(t)/ u''(t)}$ is real valued,
we see that the vector field $V$ is holomorphic if and only if
$${f'(t)\over u''(t)}=\alpha$$
for some constant $\alpha$. We note  from (2.10) that the gradient
vector field $V$ vanishes at the origin. The soliton flows along
$-V$ and hence we should expect $\alpha$ to be negative so that
everything is flowing towards the origin as in the case of the
cigar soliton when $n=1$.

It follows that
  $$ f(t)=\alpha u'(t) + k \eqno (2.11)$$ for some constant $k$.

Plugging (2.8) into (2.11), we derive the following second order
equation in $u$: $$(u')^{n-1}u''e^{\alpha u'}=\beta e^{nt} \eqno
(2.12)$$ where $\beta $ is some positive constant.

Setting $\phi (t)=u'(t)$, Eq.(2.12) becomes
$${\phi}^{n-1}{\phi}'e^{\alpha\phi}=\beta e^{nt}$$
which is a first order equation in $\phi$.

We remark that $\alpha =0$ corresponds to the flat Euclidean
metric. To get nontrivial Ricci solitons we have to require
$\alpha >0$. By scaling, we may assume $\alpha =1$. In the mean
time, we can also normalize $\beta=1$ by an appropriate
translation in $t$ (i.e., a dilation in $z$). Therefore we only
need to consider the equation
$${\phi}^{n-1}{\phi}'e^{\phi}=e^{nt} \eqno (2.13)$$

The separation of the variables $\phi$ and $t$ yields
$${\phi}^{n-1}e^{\phi}d{\phi}=e^{nt}dt$$

Integrating both sides, we get
$$\sum_{k=0}^{n-1} (-1)^{n-k-1}{n!\over k!}{\phi}^ke^{\phi}=e^{nt}+C$$

In order that $\phi\rightarrow 0$ as $t\rightarrow -\infty$ we
must have
$$C=(-1)^{n-1}n!$$ So the solution $\phi (t)$ is implicitly given by
$$\sum_{k=0}^{n-1} (-1)^{n-k-1}{n!\over k!}{\phi}^ke^{\phi}
=e^{nt}+(-1)^{n-1}n! \eqno (2.14)$$

One can check from (2.14) that the solution $\phi$   satisfies the
required asymptotic condition (2.3):
$$\phi(t)=a_1e^t+a_2e^{2t}+\cdots, \qquad a_1=1, \eqno (2.15)$$
 and the differential inequalities (2.4): $$ \phi(t)>0,\quad \phi'(t)>0,
\qquad {\hbox {for}}\ t\in (-\infty ,\infty) \eqno (2.16)$$ Thus
it gives rise to a $U(n)$-invariant K\"ahler metric $g$ on $\bf
C^n$, which is a Ricci soliton. From (2.14), or Eq.(2.13),  it is
also easy to see that
$$\lim_{t\to \infty}t^{-1}\phi(t)=n, \qquad \lim_{t\to \infty}\phi'(t)=n
\eqno (2.17)$$

To see that the soliton metric $g$ is complete, let $\rho$ denote
the distance function from the origin with respect to $g$. Since
the metric $g$ is rotationally symmetric, it is clear that
straight lines through the origin are geodesics and $\rho$ is a
function of $t$ only and is given by
$$\rho (t) =\int_{-\infty}^{t} \sqrt {\phi'(\tau) }d\tau $$
It then follows from (2.15)-(2.17) that $\rho =O(t)$, which
implies that the metric $g$ is complete.

In summary, we have proved the following

\proclaim Proposition 2.1. For each $n\ge 1$, there exists on
${\bf C^n}$ a complete rotationally symmetric K\"ahler-Ricci
soliton. Moreover, such a soliton is unique up to scaling (in the
metric) and dilation (in the variable $z$).

 Note that when $n=1$, we have $\phi =\log (1+e^t)$.
It then follows from (2.5) that
$$ds^2={|dz|^2\over 1+e^t}= {|dz|^2\over 1+|z|^2},$$
which is the cigar soliton observed by Hamilton.
\medskip

Next we shall prove that  the  soliton metric obtained above has
positive sectional curvature.

In general, the curvature tensor of a K\"ahler metric $g_{i\bar
j}$ is given by
$$R_{i\bar j k\bar l}=
-{{\partial }^2g_{i\bar j}\over \partial z^k \partial z^{\bar l}}
+g^{p\bar q}{\partial g_{i\bar q}\over \partial z^k} {\partial
g_{p\bar j}\over \partial z^{\bar l}}.$$

For a rotationally symmetric metric on ${\bf C^n}$, it suffices to
compute  the curvature at a point $P=(z_1,0,\cdots,0)$. From  a
straightforward computation we get the following:

 At point $P=(z_1,0,\cdots,0)$ ($z_1\neq 0$),
$$\eqalign { R_{i\bar j k\bar l}=e^{-2t}\Bigl \{&(\phi-\phi')({\delta}_{ij}{\delta}_{kl}
+{\delta}_{il}{\delta}_{jk})\cr &
+(3\phi'-2\phi-\phi'')({\delta}_{ij}{\delta}_{kl1}+{\delta}_{il}{\delta}_{jk1}+
{\delta}_{jk}{\delta}_{il1}+{\delta}_{kl}{\delta}_{ij1})\cr &
+(6\phi''-11\phi'-\phi'''+6\phi){\delta}_{ijkl1}
+{(\phi'-\phi'')^2\over \phi'}{\delta}_{ijkl1}\cr & +
{(\phi-\phi')^2\over \phi}({\delta}_{ij\hat 1}{\delta}_{kl1}
+{\delta}_{il\hat 1}{\delta}_{jk1}+ {\delta}_{jk\hat
1}{\delta}_{il1}+{\delta}_{kl\hat 1}{\delta}_{ij1}) \Bigr \}\cr }
\eqno (2.18)$$ where ${\delta}_{ij1}$ and ${\delta}_{ijkl1}$ are
zero unless all the indices are 1, while ${\delta}_{ij\hat 1}$ is
zero unless $i=j$ and neither index is 1.

 At the origin, the above formula reduces to
$$R_{i\bar j k\bar l}=-a_2({\delta}_{ij}{\delta}_{kl}
+{\delta}_{il}{\delta}_{jk})\eqno (2.19)$$ where $a_2$ is the
coefficient of $e^{2t}$ in (2.3).

Recall that the sectional curvature of the $2$-plane spanned by
two (real) tangent vectors
$$X= Re\ v^i{\partial \over \partial z_i} \qquad \hbox {and } \qquad Y= Re\ w^i{\partial \over \partial z_i}$$ is given by
$$||X\wedge Y||^{-2} R_{i\bar j k\bar l}(v^iw^{\bar j}-w^iv^{\bar j})
(w^kv^{\bar l}-v^kw^{\bar l})$$ where, up to a constant factor,
$$||X\wedge Y||^2=g_{i\bar l}g_{k\bar j}\bigl [(v^iw^{\bar j}-w^iv^{\bar j})
(w^kv^{\bar l}-v^kw^{\bar l})+(v^iw^k-w^iv^k) (v^{\bar l}w^{\bar
j}-w^{\bar l}v^{\bar j})\bigr ]$$ denotes the square of the area
of the parallelogram spanned by $X$ and $Y$. We remark that
$||X\wedge Y||^2=0$ (i.e., $X$ and $Y$ are colinear) if and only
if
$$v^iw^{\bar j}-w^iv^{\bar j}=0 \qquad {\hbox {for all $i$, $j$,}} \eqno (2.20)$$
Hence the positivity of sectional curvature is equivalent to
$$R_{i\bar j k\bar l}(v^iw^{\bar j}-w^iv^{\bar j})
(w^kv^{\bar l}-v^kw^{\bar l})>0$$ for all pairs of complex numbers
$(v^i)$,
 $(w^i)$
which do not satisfy (2.20). Without loss of generality, we can
assume that
$$v^i=0, \ i\geq 2 \qquad \hbox{and} \qquad w^i=0, \ i\geq 3 \eqno (2.21)$$

On one hand we have
$$R_{i\bar j k\bar l}(v^iw^{\bar j}-w^iv^{\bar j})
(w^kv^{\bar l}-v^kw^{\bar l})= 2R_{i\bar j k\bar l}v^iv^{\bar
j}w^kw^{\bar l} -2 Re R_{i\bar j k\bar l}(v^iw^{\bar j}v^kw^{\bar
l})$$

On the other hand, one computes from (2.18) and (2.21) the
holomorphic bisectional curvature
$$ R_{i\bar j k\bar l}v^iv^{\bar j}w^kw^{\bar l}=e^{-2t}\Bigl \{
\bigl [{(\phi'')^2\over \phi'}-\phi'''\bigr ]|v^1|^2|w^1|^2
+\bigl[{(\phi')^2\over \phi}-\phi''\bigr ]|v^1w^2|^2 \Bigr \}$$
and
$$R_{i\bar j k\bar l}(v^iw^{\bar j}v^kw^{\bar l})=
e^{-2t} \bigl [{(\phi'')^2\over \phi'}-\phi'''\bigr ](v^1 w^{\bar
1})^2 $$

Thus
$$\eqalign {R_{i\bar j k\bar l}&(v^iw^{\bar j}-w^iv^{\bar j})
(w^kv^{\bar l}-v^kw^{\bar l})\cr &=e^{-2t}\Bigl \{ \bigl
[{(\phi'')^2\over \phi'}-\phi'''\bigr ]|v^1w^{\bar 1}-w^1v^{\bar
1}|^2 +\bigl[{(\phi')^2\over \phi}-\phi''\bigr ]|v^1w^2|^2\Bigr
\}\cr}$$

Therefore the positivity of sectional curvature follows from the
following

\proclaim Lemma 2.2. Let $\phi $ be the solution of Eq.(2.13),
given by (2.14).
 Then we have, for all $t\in (-\infty, \infty)$,
$${\hbox {(i)}}\ \ \ \qquad \phi -{\phi}'>0$$
$${\hbox{(ii)}}\qquad ({\phi}')^2 -\phi{\phi}''>0$$
$${\hbox{(iii)}}\qquad ({\phi}'')^2 -{\phi}'{\phi}'''>0.$$

\noindent Proof.   First we note that (iii) $\Rightarrow $ (ii)
$\Rightarrow$ (i). For example, (ii) implies that the function
${\phi'/ \phi}$ is a strictly decreasing function of $t$ hence
$(i)$ follows, since $\phi'/\phi =1$ at $t=-\infty$.

We know that $\phi$ satisfies the equation
$${\phi}'={e^{nt}\over {\phi}^{n-1}e^{\phi}}\eqno (2.22)$$
On the other hand, since ${\phi' (t)}>0$ on $(-\infty, \infty)$,
$t$ can be considered as a function of $\phi$ and we obtain
$${d\over d\phi}(e^{nt})=ne^{nt}{dt\over d{\phi}}=n{\phi}^{n-1}e^{\phi}
\eqno (2.23)$$
\medskip

From (2.22) we compute that
$${\phi}''=n{\phi}'-({\phi}')^2-(n-1){({\phi}')^2\over \phi}$$ and
$$\eqalign {{\phi}'''
=&n^2{\phi}'-3n({\phi}')^2+2({\phi}')^3-3n(n-1){({\phi}')^2 \over
\phi}\cr &+4(n-1){({\phi}')^3 \over \phi}+(n-1)(2n-1){({\phi}')^2
\over {\phi }^2} \cr}$$ Hence
$$\eqalign {({\phi}'')^2 -{\phi}'{\phi}'''=&(\phi') ^2
\bigl [n\phi '-(\phi')^2+n(n-1){\phi'\over \phi}
-2(n-1){\phi'^2\over \phi}-n(n-1){(\phi')^2\over {\phi}^2}\bigr
]\cr =&{(\phi')^4\over e^{nt}{\phi}^2}\bigl [
n{\phi}^{n+1}e^{\phi}
-{\phi}^2e^{nt}+n(n-1){\phi}^ne^{\phi}-2(n-1){\phi}e^{nt}
-n(n-1)e^{nt}\bigr ]\cr}$$

Let
$$c(\phi)= n{\phi}^{n+1}e^{\phi}
+n(n-1){\phi}^ne^{\phi}-{\phi}^2e^{nt}-2(n-1){\phi}e^{nt}-n(n-1)e^{nt}
\eqno (2.24)$$ Then we have $c(0)=0$. Differentiating (2.24) with
respect to $\phi$ and using (2.23), we get
$$c'(0)=c''(0)=0$$
and
$$ c'''(\phi)=2n(n-1){\phi}^{n-2}e^{\phi}$$
Since $\phi >0$ on $(-\infty, \infty)$, this implies that
$c'''(\phi)>0$ for all $\phi >0$. This in turn implies $c''(\phi)
>0$, $c'(\phi)>0$ and $c(\phi)>0$. Thus (iii) is proved.

Moreover it follows from Lemma 2.2 (i) that the coefficient $a_2 <
0$,
 hence (2.19) implies that
the sectional curvature at the origin is also positive. Therefore,
we have proved that the gradient soliton $g$ on ${\bf C}^n$  in
Proposition 2.1 has positive sectional curvature. This finishes
the proof of Theorem 1.

\vskip 0.3in

\noindent {\bf \S 3. The Gradient Soliton on the anticanonical
bundle of ${\bf P}^{n-1}$}
\bigskip
In this section we shall consider the total space of the
anticanonical (or canonical) bundle over the complex projective
space ${\bf P}^{n-1}$ and construct a gradient K\"ahler-Ricci
soliton  on it. The soliton metric is again a rotationally
symmetric one and has nonnegative sectional curvature.

Let $z_1, z_2, \cdots, z_n$ be the coordinate system on complex
Euclidean space  ${\bf C}^n$. The complex projective space  ${\bf
P}^{n-1}$ is the quotient of ${\bf C}^n\setminus\{0\} $ by the
action of the group of nonzero scalar multiplications.  It is
covered by $n$ coordinate domains $U_j$, $j=1,2,\cdots, n$, each
characterized by $z_j\neq 0$ with the local (inhomogeneous)
coordinate system $({z_i/z_j}) (1\leq i\leq n, i\neq j)$. Let $L$
denote the hyperplane bundle over ${\bf P}^{n-1}$. For each
nonzero integer $m$, we consider the total space $X_m$ of the line
bundle $L^m \buildrel \pi \over \longrightarrow {\bf P^{n-1}}$,
where $L^m$ is the $m$th tensor power  of $L$. Note that $X_n$ is
the canonical bundle of  ${\bf P}^{n-1}$.  Without loss of
generality, we may choose $m>0$.  The transition functions of the
bundle are given by $$ y_i=\big ({z_i\over z_j}\big )^my_j \eqno
(3.1)$$ in $\pi^{-1}(U_i\bigcap U_j)$, where $y_i\in {\bf C}^1$ is
the fiber coordinate in $\pi^{-1}(U_i)$. Let $S_0$  denote the
zero cross section in $X_m$,  given by $y_i=0$. The complement
$X'_m=X_m\setminus(S_0)$ can be globally parametrized by the
homogeneous coordinate space   ${\bf C}^n\setminus\{0\} $ under
the $m$-to-one map $$(z_1, z_2, \cdots, z_n)\longrightarrow\bigl
(\bigl ({z_i\over z_j}\bigr ); (z_j)^m\bigr )\in X'_m\bigcap
\pi^{-1}(U_j) \qquad (z_j\neq 0,1\leq i\leq n; i\neq j)\eqno
(3.2)$$

As in the last section,  a  $U(n)$-invariant K\"ahler metric $g$
on  ${\bf C}^n\setminus\{O\} $ corresponds to a K\"ahler potential
function $\Phi (z,\bar z)$ given by  $$ \Phi (z,\bar z)
=u(t),\qquad \bigl(t=\log |z|^2\bigr) $$
 where $u(t)$ is a smooth function  on $(-\infty, \infty)$ which satisfies the differential inequalities (2.4). The property that the metric $g$, pulled back to  $X'_m=X_m\setminus(S_0)$, extends to a K\"ahler metric on all of $X_m$ can be translated into the following asymptotic condition on $u(t)$:

   There exists a constant $a>0$ such that  the function $u(t)-at$ has the expansion $$ u(t)-at =
a_0+a_1e^{mt}+a_2e^{2mt}+\cdots \, \eqno (3.3)$$ as $t\rightarrow
-\infty$, with $a_1>0$.

Geometrically, the positive constant $a$ times $2\pi$ will be the
area (in the metric $g$) of the complex projective line in the
zero section $S_0={\bf P}^{n-1}$.

The same argument as in section 2 (see (2.5)-(2.13)) implies that
the soliton metric must satisfy the equation
$${\phi}^{n-1}e^{\phi}d{\phi}=e^{nt}dt \eqno (3.4)$$
Hence
$$\sum_{j=0}^{n-1} (-1)^{n-j-1}{n!\over j!}{\phi}^je^{\phi}=e^{nt}+C\eqno (3.5)$$

In order that $\phi\rightarrow a$ as $t\rightarrow -\infty$ we
must have
$$ C=\sum_{j=0}^{n-1}(-1)^{n-j-1}{n!\over j!}a^je^a. \eqno (3.6)$$
Eq.(3.5) also implies that the exponent $m$ in Eq.(3.3) must be
equal to $n$, hence  the manifold $X_m$ under consideration has to
be the total space of the canonical or anticanonical line bundle
of ${\bf P}^{n-1}$.

Again from Eq.(3.4), we have

$$\lim_{t\to \infty}t^{-1}\phi(t)=n, \qquad \lim_{t\to \infty}\phi'(t)=n
\eqno (3.7)$$ which implies the soliton metric on $X_n$ is
complete.

The computation of sectional curvature in section 2 shows that
away from the zero section $S_0$ the  sectional curvature of the
soliton metric on $X_n$ is positive. 

Note that when $n=2$ the manifold $X_2$ is the tangent bundle over
the sphere ${\bf P}^1$. In this case  one can check that the
metric on each fiber has positive Gaussian curvature and is
quasi-isometric to a cigar soliton on ${\bf C}^1={\bf R}^2$. Thus
the proof of Theorem 2 is completed.

\vskip 0.3in \noindent {\bf \S 4. Homothetically Shrinking
Gradient Solitons}
\bigskip
In this section we shall use a similar argument to  construct
(homothetically shrinking) gradient K\"ahler-Ricci solitons of
Eq.(1.5) on some compact complex manifolds. These manifolds are
ones on which Calabi [2] constructed extremal metrics. They are
the total spaces of certain projective line bundles over the
complex projective space ${\bf P}^{n-1}$. Indeed we have followed
closely the notations used in [2] in the previous two sections.
Similar constructions can be carried out when the base manifold is
any compact symmetric K\"ahler manifold.

We shall use the same notations as in section 3. Let $L$ again
denote the hyperplane bundle over the projective space ${\bf
P}^{n-1}$. For each nonzero integer $k$, we consider  the total
space $M_k$ of the projective line bundle ${\bf P}(L^k\bigoplus
L^{-k})\buildrel \pi \over \longrightarrow {\bf P^{n-1}}$. The
transition functions of the bundle are given by $$ y_i=\big
({z_i\over z_j}\big )^ky_j$$ in $\pi^{-1}(U_i\bigcap U_j)$, where
$y_i\in {\bf P^1}$ is the fiber coordinate in $\pi^{-1}(U_i)$. Let
$S_0$ and $S_{\infty}$ denote the zero and $\infty$ cross sections
in $M_k$, given by $y_i=0$ and $y_i=\infty$, respectively. Without
loss of generality, we again choose $k$ to be positive. The
complement $M'_k=M_k\setminus(S_0\bigcup S_{\infty})$ can  be
globally parametrized by the homogeneous coordinate space   ${\bf
C}^n\setminus\{0\} $ under the $k$-to-one map (3.2).

From [2] we know that the maximal compact subgroup $G_K$ of the
automorphisms of $M_k$ is isomorphic to $U(n)/Z_k$. A K\"ahler
metric $g$ on $M_k$ invariant under the actions of $G_K$
corresponds to a K\"ahler metric, again denoted as $g$, on ${\bf
C}^n\setminus\{0\} $ generated by
 a K\"ahler potential   $$ \Phi (z,\bar z) =u(t),\qquad \bigl(t=\log |z|^2\bigr) \eqno (4.1)$$
 where $u(t)$ is a smooth function  on $(-\infty, \infty)$ which satisfies the differential inequalities (2.4)
and the following asymptotic properties:

 i) \  There exists a constant $a>0$ such that  the function $u(t)-at$ has the expansion $$ u(t)-at =
a_0+a_1e^{kt}+a_2e^{2kt}+\cdots \, \eqno (4.2)$$ as $t\rightarrow
-\infty$, with $a_1>0$;

ii) \  There exists a constant $b>0$ such that  the function
$u(t)-bt$ has the expansion $$ u(t)-bt =
b_0+b_1e^{-kt}+b_2e^{-2kt}+\cdots \, \eqno (4.3)$$ as
$t\rightarrow  \infty$, with $b_1<0$.

The positive constants $a, b$ actually specify the K\"ahler class
of the resulting metric. They represent the areas of the two
projective lines lying one in each of the two cross sections $S_0$
and $S_{\infty}$. We remark that the K\"ahler class of the metric
$g$ is equal to the first Chern class $C_1(M_k)$ of $M_k$ provided
$a=n-k$ and $b=n+k $. So $C_1(M_k)$ is positive precisely when
$1\leq k\leq n-1$.

For any smooth function $u(t)$, $-\infty < t<\infty$, satisfying
(2.4) and the asymptotic conditions (4.2) and (4.3) with $a=n-k$
and $b=n+k$, the formulas (2.5)-(2.10) imply that $$R_{i\bar
j}-g_{i\bar j}={\partial}_i{\partial}_{\bar j}(f-u)$$ where
$$f(t)=-\log \det (g_{i\bar j})=nt-(n-1)\log u'(t)-\log u''(t) \eqno (4.4)$$

Hence
 the gradient vector field of the function $f-u$  is given by
$$V^i=g^{i\bar j}e^{-t}z_j(f'-u')=z_i{f'(t)-u'(t)\over u''(t)}$$
so the vector field $V$ is holomorphic if and only if  $$
f'(t)-u'(t)=c_1 u''(t) \eqno (4.5)$$ for some constant $c_1$.

Plugging (4.4) into (4.5) and setting $ \phi(t)=u'(t)$, we derive
the following second order equation in $\phi$:
$${\phi''\over \phi'}+\bigl [{n-1\over \phi}+ c_1\bigr]\phi'=n-\phi \eqno (4.6)$$

Since the variable $t$ does not appear in Eq.(4.6) we can solve
for $\phi'$ and get $$\eqalign { \phi'=&{1\over
\phi^{n-1}e^{c_1\phi}}\bigl [n\int \phi^{n-1}e^{c_1\phi} d\phi -
\int \phi^ne^{c_1\phi} d\phi \bigr ]\cr =& {-1\over
c_1^{n+1}\phi^{n-1}}\Bigl [c_1^n\phi^n +
\sum_{j=0}^{n-1}(-1)^{n-j}{n!\over
j!}(1+c_1)c_1^j{\phi}^j-c_2e^{-c_1\phi}
 \Bigr ]\cr} $$ where $c_2$ is another constant.

An implicit solution  $u(t)$ with $\phi(t)=u'(t)$ is given by
$$t=-\int {c_1^{n+1}\phi^{n-1} d\phi \over c_1^n\phi^n + \sum_{j=0}^{n-1}(-1)^{n-j}{n!\over j!}(1+c_1)c_1^j{\phi}^j-c_2e^{-c_1\phi}}\eqno (4.7)$$
The asymptotic conditions (4.2) and (4.3) require that the
integrand in (4.7) has simple poles at $\phi=n-k$ and $\phi=n+k$
with residues equal to $1/k$ and $-1/k$ respectively. To have a
simple pole at $\phi=n-k$ with residue $1/k$, the constants $c_1$
and $c_2$ have to satisfy a system of two equations. It turns out
that these two equations are identical and  given by
$$ \sum_{j=0}^n(-1)^{n-j}{n!\over j!}(n-k)^{j-1}(n-k-j)c_1^j=c_2e^{-(n-k)c_1} \eqno (4.8)$$
Similarly, the condition of a simple pole at $\phi=n+k$ with
residue $-1/k$ corresponds to another  equation in $c_1$ and
$c_2$:
 $$\sum_{j=0}^n(-1)^{n-j}{n!\over j!}(n+k)^{j-1}(n+k-j)c_1^j=c_2e^{-(n+k)c_1} \eqno (4.9) $$

Hence, Eq.(4.6) admits a solution, given by (4.7), which
satisfies the conditions (4.2) and (4.3)  if and only if $c_1$ and
$c_2$ satisfy both Eq.(4.8) and Eq.(4.9). The latter condition is
easily seen to be  equivalent to that $c_1$ is a nonzero root of
the equation  $h(x)=0$, where the function $h(x)$ is given by
$$\eqalign {h(x)=&e^{2kx}\Bigl\{\sum_{j=0}^n(-1)^{n-j}{n!\over j!}(n+k)^{j-1}(n+k-j)x^j\Bigr \}\cr
&- \sum_{j=0}^n(-1)^{n-j}{n!\over j!}(n-k)^{j-1}(n-k-j)x^j. \cr }
$$

Now the existence and uniqueness of the rotationally symmetric
Ricci soliton on the manifold $M_k (1\leq k\leq n-1)$ follows from
the following

\proclaim Lemma 4.1. The equation $h(x)=0$  has one and only one
nonzero root $c_1$ with $-1<c_1<0$.

We shall only outline the proof of Lemma 4.1 and omit many
details. It turns out the proof of the uniqueness is kind of
interesting.

\smallskip
\noindent {\bf Step 1}: $h^{(i)}(0)=0$, for $0\leq i\leq n$.

From direct computations  we have, for $0\leq i\leq n$,
$$\eqalign { h^{(i)}(0)=&\sum_{j=0}^i
{i\choose j}(2k)^{i-j}(-1)^{n-j}n!(n+k)^{j-1} (n+k-j)\cr
&-(-1)^{n-i}n!(n-k)^{i-1}(n-k-i)\cr =&(-1)^{n-i}n! \bigl
[(n-k)^i-i(n-k)^{i-1}-(n-k)^{i-1}(n-k-i)\bigr ]=0\cr}$$
\smallskip
 \noindent {\bf Step 2}: $h^{(n+1)}(x)>0$, for $x\geq 0$.

From direct computations we have

$$ \eqalign {h^{(n+1)}(x)=&e^{2kx}\sum_{j=0}^n {n+1\choose j}(2k)^{n+1-j}\sum_{i=0}^{n-j} (-1)^{n-i-j}{n!\over i!}(n+k)^{i+j-1}(n+k-i-j)x^i\cr
=&e^{2kx}\sum_{i=0}^n C_ix^i\cr}$$ where the coefficients $C_i,
i=0,1, ..., n$, are given by

$$C_i={n!\over i!}(n+k)^i\Bigl [\sum_{j=0}^{n-i}{n+1\choose j}(-1)^{n-i-j}(2k)^{n+1-j}(n+k)^{j-1} (n+k-i-j)\Bigr ].$$

We claim that $C_i>0$ for $0\leq i\leq n$. To see this, let
$$B_j={n+1\choose j}(2k)^{n+1-j}(n+k)^{j-1}(n+k-i-j).$$
Then it is easy to show that $$B_{j+1}>B_j, \qquad {\hbox {for}}\
0\leq j\leq  n-i-1\eqno (4.11)$$ Note that for each $i$, $C_i$ is
an alternating sum starting with a positive term when $j=n-i$.
Therefore (4.11) implies that $C_i>0$.
\smallskip
 \noindent {\bf Step 3}: The sign of $h(-1)$ is $(-1)^n$.

From direct computations we have
$$\eqalign {h(-1)=&(-1)^n n!\Bigl [e^{-2k}\sum_{j=0}^n{1\over j!}(n+k)^{j-1}(n+k-j)-
\sum_{j=0}^n{1\over j!}(n-k)^{j-1}(n-k-j)\Bigr ]\cr =& (-1)^n
\bigl [e^{-2k}(n+k)^n-(n-k)^n\bigr ].\cr}$$ Now it can be shown
that $$e^{-2k}(n+k)^n>(n-k)^n, \qquad {\hbox {for}}\ 1\leq k\leq
n-1.$$

It follows from Step 1--Step 3 that $h$ has no zero in $(0,
\infty)$ and has a zero in the interval $(-1, 0)$.

\smallskip
 \noindent {\bf Step 4}: $h$ has only one zero in the interval $(-\infty, 0)$.

Let $y=-x$, then $h(x)=(-1)^{n+1}e^{-2ky} g(y)$ with
$$ g(y)=e^{2ky}\Bigl\{\sum_{j=0}^n{n!\over j!}(n-k)^{j-1}(n-k-j)y^j\Bigr \}
- \sum_{j=0}^n{n!\over j!}(n+k)^{j-1}(n+k-j)y^j.  \eqno (4.12) $$

It is then equivalent to show that $g$  only has one zero in $(0,
\infty)$. The crucial point here is to observe that on the
interval $(0, \infty)$, $g$ can be written as a power series of
the following special form:
$$g(y)=\sum_{j=n+1}^{n+l}a_jy^j- \sum_{j=n+l+1}^{\infty}a_jy^j \eqno (4.13)$$ where $a_j>0$ for all $j$ and $l\geq 2$ is some positive integer. Note that for any function $g$ of the type (4.13), the number of sign changes of $g$ in the interval $(0, \infty)$ must agree with the number of sign changes of its derivatives $g^{(k)}, 1\leq k\leq n+l$. Hence it can have only one zero in  $(0,\infty)$.

To see that $g$ has the form in (4.13), we compute the Taylor
expansion of $g$:
$$g(y)=n!\sum_{i=n+1}^{\infty}b_iy^i=n!\sum_{i=n+1}^{\infty} \Bigl
\{\sum_{j=0}^n{(2k)^{i-j}(n-k)^{j-1}\over (i-j)! j!}(n-k-j)\Bigr
\}y^i.
$$

It can be shown that  $\{b_i\}$ is a decreasing sequence. Also
$b_1$ must be positive yet $b_i$ cannot all be positive because
$g(y)$ does have a zero in $(0, \infty)$.  Hence the Taylor
expansion of $g$ must have the form stated in (4.13). Thus the
proof of Lemma 4.1 is completed.

\medskip

Finally we  mention that the Ricci curvature of the soliton metric
$g$ on $X_k$ is positive if and only if  $k=1$. From (2.9) we have
$$R_{i\bar j}=e^{-t}f'(t){\delta}_{i\bar j}
+e^{-2t}{\bar z}_iz_j\bigl(f''(t)-f'(t)\bigr). \eqno (4.10)$$
Hence $R_{i\bar j}>0$ if and only if  $f'>0$ and $f''>0$, or by
(4.5), if and only if $\phi +c_1\phi' >0$ and $\phi'
+c_1\phi''>0$. While the first inequality is true for $1\leq k\leq
n-1$, the second one holds only when $k=1$.

\medskip

\noindent {\bf Notes added for the arXiv posting:} This is the original paper appeared in the book ``Elliptic and Parabolic Methods in Geometry (Minneapolis, MN,
1994), A K Peters, Wellesley, MA, (1996)"  (p.1-16), except with the following modifications: 

\smallskip

1.  Corrected an inaccuracy in the statement of Theorem 2 by removing the phrase 

\hskip 0.16 in   ``non-negativity  of sectional  curvature" (of the soliton metric); 

2.  Modified Remark 3 and added references [9] and [10];
 
3.  Added Theorem 3 in the introduction to summarize the main result in Section 4.

\beginsection References \par
\medskip
\frenchspacing

\item{[1]} Bryant, R.: Existence of a gradient Ricci soliton in
dimension three. Preprint.
\smallskip

\item{[2]} Calabi, E.: Extremal K\"ahler metrics. Seminars on
Differentail Geometry (S.-T. Yau ed.), Princeton Univ. Press \&
Univ. of Tokyo Press, Princeton, New York, 1982, pp. 259-290
\smallskip
\item{[3]} Cao, H.-D.: On Harnack's inequalities for the
K\"ahler-Ricci flow.  Invent. Math. {\bf 109}, 1993, pp. 247-263
\smallskip

\item{[4]} Hamilton, R.S.: The Ricci flow on surfaces. in
Mathematics and General Relativity, Contemporary Mathematics {\bf
71}, 1986.
\smallskip

\item{[5]} Hamilton, R.S.: The Harnack estimate for the Ricci
flow. J. Differ.Geom. {\bf 37}, 1993, pp. 225-243
\smallskip

\item{[6]} Hamilton, R.S.: Eternal solutions to the Ricci flow.
Preprint

\item{[7]} Ivey, T.: Ricci solitons on Compact three-manifolds.
Differ. Geom. and its Appl. {\bf 3}, 1993, pp. 301-307

\item {[8]} Koiso, N: On rotationally symmetric Hamilton's
equation for K\"ahler-Einstein metrics. Advanced Studies in Pure
Mathematics 18-I, 1990, Recent Topics in Differential and Analytic
Geometry, pp. 327-337

\item{[9]} Calabi, E.:   M\'etriques k\"ahl\'eriennes et fibr\'es
olomorphes,  Ann. Sci. \'Ecole Norm. Sup., {\bf 12} (1979),
269--294

\item{[10]} Eguchi, T. and  Hanson, A. J.:   Asymptotically flat
selfdual solutions to Euclidean gravity, Phys. Lett. B
{\bf 74} (1978), 249--251

\vskip 0.3in \noindent Department of Mathematics
\smallskip
\noindent Texas A\&M University
\smallskip
\noindent College Station, TX 77843
\smallskip \noindent Email: cao@math.tamu.edu

\end